\documentclass[12pt]{article}
\usepackage{latexsym,amsmath,amsthm}
\usepackage{amssymb}
\usepackage{url}
\usepackage[margin=1in]{geometry}
\usepackage{tikz}

\newtheorem{thm}{Theorem}[section]
\newtheorem{prop}[thm]{Proposition}
\newtheorem{question}[thm]{Question}
\newtheorem{cor}[thm]{Corollary}

\theoremstyle{definition}
\newtheorem{prob}[thm]{Problem}

\theoremstyle{definition}
\newtheorem{defn}[thm]{Definition}

\theoremstyle{remark}
\newtheorem{remark}[thm]{Remark}

\newcommand{\currentbestpath}{5000}

\tikzstyle{W}=[draw,circle, fill=white, minimum size=6pt,inner sep=0pt]

\begin{document}

\title{The independent set sequence of some families of trees}

\author{David Galvin\thanks{Department of Mathematics, 255 Hurley Hall, University of Notre Dame,
Notre Dame IN 46556, USA; dgalvin1@nd.edu. Research supported by NSA grant H98230-13-1-0248, and by the Simons Foundation.}~~and Justin Hilyard\thanks{Department of Mathematics, 255 Hurley Hall, University of Notre Dame,
Notre Dame IN 46556, USA;  jhilyard@alumni.nd.edu. Research supported by NSA grant H98230-13-1-0248.}}

\date{\today}

\maketitle

\begin{abstract}
For a tree $T$, let $i_T(t)$ be the number of independent sets of size $t$ in $T$. It is an open question, raised by Alavi, Malde, Schwenk and Erd\H{o}s, whether the sequence $(i_T(t))_{t \geq 0}$ is always unimodal. Here we answer the question in the affirmative for some recursively defined families of trees, specifically paths with auxiliary trees attached at the vertices in a periodic manner. In particular, extending a result of Wang and B.-X. Zhu, we show unimodality of the independent set sequence of a path on $2n$ vertices with $\ell_1$ and $\ell_2$ pendant edges attached alternately at the vertices of the path, $\ell_1, \ell_2$ arbitrary. 

We also show that the independent set sequence of any tree becomes unimodal if sufficiently many pendant edges are attached at any single vertex, or if $k$ pendant edges are attached at every vertex, for sufficiently large $k$. This in particular implies the unimodality of the independent set sequence of some non-periodic caterpillars. 
 
\end{abstract}

\section{Introduction and statement of results} \label{sec-intro}

An {\em independent set} (stable set) in a graph $G$ is a set of pairwise non-adjacent vertices. Denote by $i_G(t)$ the number of independent sets in $G$ of size $t$ (with $t$ vertices), and by $\alpha(G)$ the size of the largest independent set in $G$. The {\em independent set sequence} of $G$ is the sequence $(i_G(t))_{t=0}^{\alpha(G)}$. (All graphs in this note are simple, finite and loopless.)

A seminal result of Heilmann and Lieb \cite{HeilmannLieb} on the matching polynomial of a graph implies that if $G$ is a line graph then the independence polynomial of $G$ has the real-roots property and so the independent set sequence is log-concave and unimodal (see Definition \ref{defn-basic_terms} and Remark \ref{remark-Newton} below). For the more general class of claw-free graphs (graphs without an induced star on four vertices), log-concavity and unimodality of the independent set sequence was shown by Hamidoune \cite{Hamidoune}, and later Chudnovsky and Seymour demonstrated the stronger real-roots property for this family \cite{ChudnovskySeymour}.

\begin{defn} \label{defn-basic_terms}
A finite sequence $(a_0, a_1, \ldots, a_n)$ of real numbers has the {\em real-roots property} if the generating function of the sequence (the polynomial $\sum_{k=0}^n a_kx^k$) factors into $n$ linear terms over the reals. The sequence is {\em log-concave} if $a_k^2 \geq a_{k-1}a_{k+1}$ for $k=1, \ldots, n-1$, and is {\em unimodal} if there is some $m$, $0 \leq m \leq n$, such that $a_0 \leq a_1 \leq \ldots \leq a_m \geq a_{m+1} \geq \ldots \geq a_n$. We also say that a generating function of a finite sequence is log-concave (or unimodal), if its coefficient sequence is log-concave (or unimodal). 
\end{defn}
We refer to the generating function of the independent set sequence of a graph $G$, an object introduced by Gutman and Harary \cite{GutmanHarary}, as the {\em independence polynomial} of $G$, and denote it by $p(G,x)$.    

\begin{remark} \label{remark-Newton}
Let $(a_0, a_1, \ldots, a_n)$ be a sequence of positive numbers. If it has the real roots property then it is log-concave; see e.g. \cite[Chapter 8]{Bona}. If it is log-concave, in which case we say that both it and its generating function are {\em LC$^+$}, then it is unimodal.
\end{remark}

Alavi, Malde, Schwenk and Erd\H{o}s \cite{AlaviMaldeSchwenkErdos}, considering a question of Wilf, showed that in general the independent set sequence can exhibit essentially any pattern of rises and falls. Specifically, they exhibited, for each $m \geq 1$ and each permutation $\pi$ of $\{1, \ldots, m\}$, a graph $G$ with $\alpha(G)=m$ for which
$$
i_G(\pi(1)) < i_G(\pi(2)) < \ldots < i_G(\pi(m)).
$$
They then considered the question of whether there are other families, besides claw-free graphs, with unimodal independent set sequence. Paths, being claw-free, certainly have unimodal independent set sequence, and it is easy to verify that the same is true for stars, the other natural extremal family of trees. Perhaps based on these observations Alavi et al. posed an intriguing question that is the subject of the present paper.
\begin{question} \label{quest-AEMS}
Is the independent set sequence of every tree unimodal?
\end{question}
Despite substantial effort, not much progress has been made on this question since it was raised in 1987.
We briefly review here some of the families of trees for which the unimodality of the independent set sequence has been established. 
\begin{itemize}
\item A {\em spider} is a tree with at most one vertex of degree at least 3, and a graph is {\em well-covered} if all of its maximal independent sets have the same size. Levit and Mandrescu \cite{LevitMandrescu2} showed that all well-covered spiders have unimodal independent set sequence.
\item A tree on $n$ vertices is {\em maximal} if it has the greatest number of maximal independent sets (with respect to inclusion) among $n$ vertex trees; maximal trees were characterized by Sagan \cite{Sagan}, and belong to a class of graphs known as {\em batons}. Mandrescu and Spivak \cite{MandrescuSpivak} showed that all maximal trees with an odd number of vertices have unimodal independent set sequence, and they showed the same for some maximal trees with an even number of vertices. 
\end{itemize}

Many of the families for which unimodality of the independent set sequence has been established have the following recursive structure.
\begin{defn} \label{defn-vertebrated}
Let $G$ be a graph with a distinguished vertex $v$. An {\em $n$-concatenation} of $G$ is obtained by taking $n$ vertex-disjoint copies of $G$, with, say, the distinguished vertices labeled $v_1, v_2, \ldots, v_n$, and adding the edges $v_1v_2, v_2v_3, \ldots, v_{n-1}v_n$. (This is a special case of the rooted product construction introduced by Godsil and McKay \cite{GodsilMcKay}.)
\end{defn}

\begin{itemize}
\item An {\em $n$-centipede} is a path on $n$ vertices with a pendant edge attached at each vertex (see Figure \ref{fig-centipede}), or equivalently an $n$-concatenation of $K_2$. Levit and Mandrescu \cite{LevitMandrescu} showed that all centipedes have unimodal independent set sequence, and later Z.-F. Zhu \cite{Zhu} showed that the sequence has the real-roots property in this case. Z.-F. Zhu \cite{Zhu} also considered an $n$-concatenation of the star $K_{1,2}$ with the vertex of degree $2$ taken as the distinguished vertex (i.e., the family of trees obtained from paths by attaching two pendant edges to each vertex), and showed that all these trees have unimodal independent set sequence.

\begin{figure}[ht!] \center
\begin{tikzpicture}
	\node[W] (0) {};
	\node[right of=0, W] (1) {};
    \node[right of=1, W] (2) {};
    \node[right of=2, W] (3) {};
    \node[right of=3, W] (4) {};
    \node[right of=4, W] (5) {};
    \node[right of=5, W] (6) {};
    \node[right of=6, W] (7) {};
    \node[below of=0, W] (8) {};
    \node[below of=1, W] (9) {};
    \node[below of=2, W] (10) {};
    \node[below of=3, W] (11) {};
    \node[below of=4, W] (12) {};
    \node[below of=5, W] (13) {};
    \node[below of=6, W] (14) {};
    \node[below of=7, W] (15) {}; 	
	\draw (0) edge (1);
    \draw (1) edge (2);
    \draw (2) edge (3);
    \draw (3) edge (4);
    \draw (4) edge (5);
    \draw (5) edge (6);
    \draw (6) edge (7);
    \draw (0) edge (8);
    \draw (1) edge (9);
    \draw (2) edge (10);
    \draw (3) edge (11);
    \draw (4) edge (12);
    \draw (5) edge (13);
    \draw (6) edge (14);
    \draw (7) edge (15);
\end{tikzpicture}
\caption{An $8$-centipede.}
\label{fig-centipede}
\end{figure}
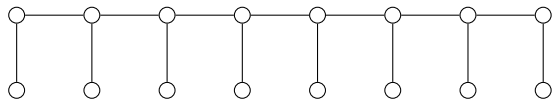

\item Wang and B.-X. Zhu \cite{WangZhu}, generalizing the work of Z.-F. Zhu and of Levit and Mandresecu, considered an $n$-concatenation of the star $K_{1,k}$, both with the vertex of degree $k$ taken as the distinguished vertex, and with a vertex of degree 1, and showed that in both cases for all $k \geq 1$ these trees have unimodal independent set sequence. 

\item B.-X. Zhu, generalizing another result from \cite{WangZhu}, obtained the following. 
\begin{prop} \cite[Theorem 3.3]{Zhu2} \label{prop-zhu}
Let $G$ be a graph (not necessarily a tree) whose independent set sequence has the real-roots property, and let $H$ be a claw-free graph (not necessarily a tree) with distinguished vertex $v$. Let $G_v[H]$ be obtained by taking $|V(G)|$ vertex-disjoint copies of $H$, with the distinguished vertices labeled $v_1, v_2, \ldots, v_{|V(G)|}$, and adding edges between the $v_i$'s so that the subgraph induced by the $v_i$'s is isomorphic to $G$. Then the independent set sequence of $G_v[H]$ has the real-roots property.
\end{prop}
The graph $G_v[H]$ is also known as the {\em rooted product} of $G$ and $H$ 
\cite{GodsilMcKay}. 
The consequence of Proposition \ref{prop-zhu} for Question \ref{quest-AEMS} is that if $T$ is any tree whose independent set sequence has the real-roots property, and if $T'$ is obtained from $T$ by attaching a path of a fixed length at each vertex (with the point of attachment being anywhere along the path, as long as the same point of attachment is chosen for each path), then the independent set sequence has the real-roots property and so is unimodal.

Starting with the family of paths and closing under the operation of attaching fixed length paths at each vertex, we obtain a large family of trees, most of which are not very ``path-like,'' with unimodal independent set sequence (the unimodality of this particular family is also established --- using different methods --- in \cite{BahlsBaileyOlsen}).

\end{itemize}

One partial result valid for all trees has been obtained. Levit and Mandrescu \cite{LevitMandrescu3} showed that if $G$ is a tree then the final one third of its independent set sequence is decreasing:
\begin{equation} \label{LevitMandrescu-finalthird}
i_{\lceil (2\alpha(G)-1)/3\rceil}(G) \geq i_{\lceil (2\alpha(G)-1)/3\rceil+1}(G) \geq \ldots \geq i_{\alpha(G)}(G).
\end{equation}

\begin{remark}
Levit and Mandrescu showed that (\ref{LevitMandrescu-finalthird}) holds for all $G$ in the class of {\em K\"onig-Egerv\'ary} graphs (in which the size of the largest independent set plus the size of the largest matching equals the number of vertices), which includes not just trees but bipartite graphs. They made the conjecture that all K\"onig-Egerv\'ary graphs have unimodal independent set sequence, but a bipartite counterexample was found by Bhattacharyya and Kahn \cite{BhattacharyyaKahn}. 
\end{remark}

\medskip

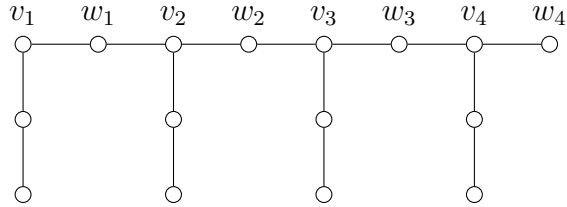
\begin{figure}[ht!] \center
\begin{tikzpicture}
	\node[W, label=above:{$v_1$}] (1) {};
    \node[right of=1, W, label=above:{$w_1$}] (2) {};
    \node[right of=2, W, label=above:{$v_2$}] (3) {};
    \node[right of=3, W, label=above:{$w_2$}] (4) {};
    \node[right of=4, W, label=above:{$v_3$}] (5) {};
    \node[right of=5, W, label=above:{$w_3$}] (6) {};
    \node[right of=6, W, label=above:{$v_4$}] (7) {};
    \node[right of=7, W, label=above:{$w_4$}] (8) {};
    \node[below of=1, W] (10) {};
    \node[below of=10, W] (11) {};
    \node[below of=3, W] (30) {};
    \node[below of=30, W] (31) {};
    \node[below of=5, W] (50) {};
    \node[below of=50, W] (51) {};
    \node[below of=7, W] (70) {};
    \node[below of=70, W] (71) {}; 	
	\draw (1) edge (2);
    \draw (2) edge (3);
    \draw (3) edge (4);
    \draw (4) edge (5);
    \draw (5) edge (6);
    \draw (6) edge (7);
    \draw (7) edge (8);
    \draw (1) edge (10);
    \draw (10) edge (11);
    \draw (3) edge (30);
    \draw (30) edge (31);
    \draw (5) edge (50);
    \draw (50) edge (51);
    \draw (7) edge (70);
    \draw (70) edge (71);
\end{tikzpicture}
\caption{A tree $T^4(v,w)$, where $T$ is a path on three vertices, with $w$ a leaf and $v$ its unique neighbor.}
\label{fig-2-centipede}
\end{figure}

The approach of Z.-F. Zhu and of Wang and B.-X. Zhu to independence polynomials has been developed considerably by Bahls and Salazar \cite{BahlsSalazar}, Bahls \cite{Bahls} and Bahls, Bailey and Olsen \cite{BahlsBaileyOlsen}, but for the most part this development does not address trees. One aim of the present note is to consider some trees that are a natural modification to the family of concatenated graphs dealt with by Levit and Mandrescu, Z.-F. Zhu, and Wang and B.-X. Zhu.
\begin{defn} \label{defn-mainobject}
Let $G$ be a graph with two distinguished vertices $v$ and $w$ that are adjacent.
An {$n$-concatenation of $G$ through $v$ and $w$}, which we denote by $G^n(v,w)$, is obtained by taking $n$ vertex-disjoint copies of $G$, with, say, the distinguished vertices labeled $v_1, w_1, v_2, w_2, \ldots, v_n, w_n$, and adding the $n-1$ edges $w_1v_2, w_2v_3, \ldots, w_{n-1}v_n$.
\end{defn}
See Figure \ref{fig-2-centipede}. Denote by $N[u]$ the closed neighborhood of a vertex $u$ ---  that is, the vertex $u$ together with the set of vertices adjacent to $u$ --- in whatever graph is under discussion, and denote by $G_a$ the graph $G-N[a]$ for any vertex $a$. 

\begin{thm} \label{thm-tech}
Let $G$ be a graph (not necessarily a tree) and let $v$ and $w$ be two adjacent vertices of $G$. Suppose that
\begin{itemize}
\item $p(G,x)$ is LC$^+$ (see Remark \ref{remark-Newton}) and 
\item $p^2(G,x) - 4qx^2p(G_v,x)p(G_w,x)$ is LC$^+$ for all $q\in[0,1]$.
\end{itemize}
Then for all $n \geq 0$ the independent set sequence of $G^n(v,w)$ (which is a tree if $G$ is) is log-concave, and hence unimodal.
\end{thm}   

Theorem \ref{thm-tech} may be used in an {\em ad hoc} manner to establish unimodality of the independent set sequence of paths of odd length that have pendant trees $T_1$, $T_2$ attached alternately, by taking $T$ to be the tree obtained from $T_1$ and $T_2$ by adding an edge joining the vertices at which these trees are attached to the paths. For example we can establish the unimodality of the independent set sequence of the $2n$-centipede ($T_1$, $T_2$ both an edge) by applying Theorem \ref{thm-tech} with $T$ a path on four vertices and with $v$ and $w$ the two non-leaf vertices.  We have
$p(T,x) = 1+4x+3x^2$ (which is evidently LC$^+$)
and
$$
p^2(T,x) - 4qx^2p(T_v,x)p(T_w,x) = 1+8x+(22-4q)x^2+(24-8q)x^3+(9-4q)x^4,
$$
which is easily seen to be LC$^+$ for all $q\in[0,1]$.

\begin{figure}[ht!] \center
\begin{tikzpicture}
	\node[W, label=above:{$v_1$}] (0) {};
	\node[right of=0, W, label=above:{$w_1$}] (1) {};
    \node[right of=1, W, label=above:{$v_2$}] (2) {};
    \node[right of=2, W, label=above:{$w_2$}] (3) {};
    \node[right of=3, W, label=above:{$v_3$}] (4) {};
    \node[right of=4, W, label=above:{$w_3$}] (5) {};
    \node[right of=5, W, label=above:{$v_4$}] (6) {};
    \node[right of=6, W, label=above:{$w_4$}] (7) {};
    \node[below of=0, W] (8) {};
    \node[below left of=1, W] (9) {};
    \node[below right of=1, W] (9a) {};
    \node[below of=2, W] (10) {};
    \node[below left of=3, W] (11) {};
    \node[below right of=3, W] (11a) {};
    \node[below of=4, W] (12) {};
    \node[below left of=5, W] (13) {};
    \node[below right of=5, W] (13a) {};
    \node[below of=6, W] (14) {};
    \node[below left of=7, W] (15) {}; 	
    \node[below right of=7, W] (15a) {};
	\draw (0) edge (1);
    \draw (1) edge (2);
    \draw (2) edge (3);
    \draw (3) edge (4);
    \draw (4) edge (5);
    \draw (5) edge (6);
    \draw (6) edge (7);
    \draw (0) edge (8);
    \draw (1) edge (9); \draw (1) edge (9a);
    \draw (2) edge (10);
    \draw (3) edge (11); \draw (3) edge (11a);
    \draw (4) edge (12);
    \draw (5) edge (13); \draw (5) edge (13a);
    \draw (6) edge (14);
    \draw (7) edge (15); \draw (7) edge (15a);
\end{tikzpicture}
\caption{A tree $S_{1,2}^4(v,w)$, where $S_{1,2}$ is a double star, $v$ is the center vertex with one pendant edge and $w$ is the center vertex with two pendant edges.}
\label{fig-mixed-star}
\end{figure}
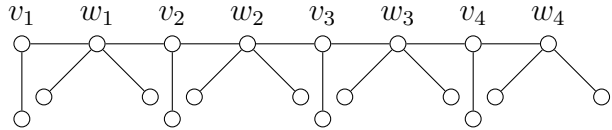

One purpose of this note is to use Theorem \ref{thm-tech} to verify the log-concavity (and hence unimodality) of some infinite families of trees that consist of paths with pairs of trees attached alternately at the vertices. Item \ref{main-thm-item-1} below subsumes a previously mentioned result of Wang and B.-X. Zhu \cite{WangZhu} (the case $\ell_1=\ell_2$).
\begin{thm} \label{maintheorem}
For each of the following trees $T$, with given distinguished vertices $v$ and $w$, the independent set sequence of $T^n(v,w)$ is log-concave (and so unimodal) for all $n \geq 0$.
\begin{enumerate}
\item \label{main-thm-item-1} $T=S_{\ell_1,\ell_2}$, a double star consisting of adjacent vertices $v$ and $w$ with $v$ having $\ell_1$ neighbors (other than $w$), all pendant edges, and $w$ having $\ell_2$ neighbors (other than $v$), all pendant edges, $\ell_1, \ell_2 \geq 0$ arbitrary. In this case $S_{\ell_1,\ell_2}^n(v,w)$ is a path on $2n$ vertices with $\ell_1$ and $\ell_2$ edges attached alternately at the vertices. See Figure \ref{fig-mixed-star}.
\item \label{main-thm-item-2} $T=P_k$ is a path on $k$ vertices, $1 \leq k \leq \currentbestpath$, with $w$ a leaf and $v$ its unique neighbor. In this case $P^n_k(v,w)$ is a path on $2n$ vertices, with a path of length $k-2$ attached at every second vertex, $k \leq \currentbestpath$. See Figure \ref{fig-2-centipede}.
\end{enumerate}  
\end{thm}

The proof of Theorem \ref{thm-tech} appears in Section \ref{sec-general_approach} and we apply it to obtain Theorem \ref{maintheorem} in Section \ref{sec-app}. Our approach follows lines similar to those of Wang and B.-X. Zhu. We obtain a recurrence for the independence polynomial $p(G^n(v,w),x)$ of $G^n(v,w)$, solve to obtain an explicit expression for $p(G^n(v,w),x)$, factorize this into bounded-degree real polynomials, and then work to establish log-concavity of each of the factors. The log-concavity of $p(G^n(v,w),x)$ then follows from the well-known fact that the product of log-concave polynomials is log-concave (see e.g. \cite[Chapter 8]{Bona}). 

We encounter an obstacle, however, not seen in the cases treated by Wang and B.-X. Zhu: there are some instances of trees $T$ covered by Theorem \ref{maintheorem} where the factors in the natural factorization of $p(T^n(v,w),x)$ are not all log-concave. Here we come to a novelty of the present work. In these cases we are still able to conclude log-concavity of the product polynomial, by clustering the factors in such a way that the product of the factors within each cluster is log-concave. 

Let us remark here that for item \ref{main-thm-item-2} of Theorem \ref{maintheorem}, for each fixed $k$ we reduce the verification that for all $n$, $P^n_k(v,w)$ has log-concave independent set sequence, to the verification that the coefficient sequence of a single polynomial of degree roughly $k$ satisfies certain properties (see the discussion around (\ref{eq-path-ver})). We have not found an argument that allows us to treat this polynomial for general $k$, so performed the verification via a {\tt Mathematica} computation. The upper bound on $k$ in Theorem \ref{maintheorem}, item \ref{main-thm-item-2} reflects the point where we stopped our computations.  

\medskip

We end the note by providing two more schemes for producing families of trees that have unimodal independent set sequence. Both involve augmenting a graph by attaching pendant stars. The first is somewhat related to Proposition \ref{prop-zhu}.
\begin{prop} \label{prop-claw}
Let $G$ be any graph and let $G_k$ be obtained from $G$ by attaching $k$ pendant edges at each vertex. For all sufficiently large $k$ (how large depending on $G$), $G_k$ has log-concave (and hence unimodal) independent set sequence.
\end{prop}  
We present the short proof in Section \ref{sec-rr}. Note that unlike Proposition \ref{prop-zhu}, here we do not require the independence polynomial of $G$ to have the real-roots property. 

Our second scheme for producing families of trees that have unimodal independent set sequence is encapsulated in the following result.
\begin{thm} \label{thm-pendant}
Let $G$ be an arbitrary graph (not necessarily a tree), and let $v$ be an arbitrary vertex of $G$. Let $G^1_n$ be obtained from $G$ by attaching a star with $n$ leaves to $G$ at $v$, with the center of the star as the point of attachment (i.e., by attaching $n$ pendant edges at $v$), and let $G^2_n$ be obtained from $G$ by attaching a star with $n$ leaves to $G$ at $v$, with one of the leaves of the star as the point of attachment (i.e., by attaching a pendant edge at $v$ to new vertex $w$, and then attaching $n-1$ pendant edges at $w$). For all sufficiently large $n=n(G)$, both $G^1_n$ and $G^2_n$ have log-concave (and hence unimodal) independent set sequence.
\end{thm}  
We give the proof in Section \ref{sec-pendant}. Here we mention a corollary.
A {\em $(k_1, \ldots, k_n)$-caterpillar} is a path on $n$ vertices with $k_i$ pendant edges attached at the $i$th vertex; so one of the previously mentioned results of Wang and B.-X. Zhu \cite{WangZhu} deals with the subfamily of $(k,k, \ldots, k)$-caterpillars, while item \ref{main-thm-item-1} of Theorem \ref{maintheorem} deals with the larger subfamily of $(\ell_1,\ell_2, \ldots, \ell_1,\ell_2)$-caterpillars. Caterpillars are a natural common extension of paths ($(0,0,\ldots,0)$-caterpillars) and stars ($(k)$-caterpillars). Note that if $G$ is a $(k_1, \ldots, k_{n-1})$-caterpillar with $v$ the vertex with $k_{n-1}$ pendant edges then $G^2_{k_n+1}$ is a $(k_1, \ldots, k_n)$-caterpillar.
\begin{cor}
For all $n$ and $k_1,\ldots, k_{n-1}$, if $k_n$ is sufficiently large then the $(k_1, \ldots, k_n)$-caterpillar has unimodal independent set sequence.   
\end{cor}

\section{Proof of Theorem \ref{thm-tech}} \label{sec-general_approach}

Recall that $N[u]$ is the closed neighborhood of a vertex $u$ in whatever graph is under discussion. We have an easy identity:
\begin{equation} \label{rec}
p(H,x) = p(H-ab,x) - x^2p(H-N[a]-N[b],x)
\end{equation}
for any graph $H$ and any edge $ab$ of $H$. 

Let $G$ be given, with distinguished vertices $v$ and $w$ that are adjacent. For typographic clarity, denote by $p_n(x)$ the independence polynomial $p(G^n(v,w),x)$ of $G^n(v,w)$, and recall that $G_a$ is the graph $G-N[a]$ for any vertex $a$. Applying (\ref{rec}) to $G^n(v,w)$ with $a=w_{n-1}$ and $b=v_n$ we obtain the recurrence relation
\begin{equation} \label{eq-main_recurrence}
p_n(x) = p(G,x)p_{n-1}(x) - x^2p(G_v,x)p(G_w,x)p_{n-2}(x),
\end{equation}
for $n \geq 2$ with initial conditions $p_0(x) = 1$ and $p_1(x) = p(G,x)$. We use here that $G^n(v,w)-w_{n-1}v_n$ consists of a copy of $G^{n-1}(v,w)$ and a copy of $G$ with no edges between them (so the independence polynomial of $G^n(v,w)$ is the product of those of $G^{n-1}(v,w)$ and $G$), and that $G^n(v,w)-N[w_{n-1}]-N[v_n]$ consists of copies of $G^{n-2}(v,w)$, $G_w$ and $G_v$ with no edges between them; that there are no edges between the copies of $G^{n-2}(v,w)$ and $G_w$ uses that $v$ and $w$ are adjacent.

We can explicitly solve this recurrence using standard methods. From \cite{Brualdi} we have that if $(z_n)_{n \geq 0}$ is a sequence satisfying $z_n=az_{n-1}+bz_{n-2}$ for $n \geq 2$ with $a^2+4b > 0$ then
$$
z_n = \frac{(z_1-z_0\lambda_2)\lambda_1^n + (z_0\lambda_1-z_1)\lambda_2^n}{\lambda_1-\lambda_2}
$$
for $n \geq 0$, where $\lambda_1 = (a+\sqrt{a^2+4b})/2$ and $\lambda_2 = (a-\sqrt{a^2+4b})/2$ are the roots of $\lambda^2-a\lambda -b =0$. Applying to the present situation, where $z_1=a=p(G,x)$, $z_0=1$ and $b=-x^2p(G_v,x)p(G_w,x)$, and using $\lambda_1+\lambda_2=a$, we obtain
$$
p_n(x) = \frac{\lambda_1^{n+1}-\lambda_2^{n+1}}{\lambda_1-\lambda_2}.
$$ 
Note that having $z_1/z_0=a$ is critical here for obtaining a clean final expression for $p_n(x)$.

What follows is similar to \cite[(3.1)]{WangZhu}, and can be derived from a combination of Lemmas 2.3 and 2.4 of \cite{WangZhu}. Using 
$$
\lambda_1^{n+1}-\lambda_2^{n+1} = \left\{
\begin{array}{ll}
(\lambda_1-\lambda_2)(\lambda_1+\lambda_2) \prod_{s=1}^{(n-1)/2} \left((\lambda_1+\lambda_2)^2-4\lambda_1\lambda_2\cos^2\left(\frac{s\pi}{n+1}\right)\right) & \mbox{if $n$ odd} \\
(\lambda_1-\lambda_2)\prod_{s=1}^{n/2} \left((\lambda_1+\lambda_2)^2-4\lambda_1\lambda_2\cos^2\left(\frac{s\pi}{n+1}\right)\right) & \mbox{if $n$ even}
\end{array}
\right.
$$
(see e.g. \cite{BarnardChild}) along with 
$\lambda_1\lambda_2=-b$, we obtain 
\begin{equation} \label{eq-explicit_expression}
p_n(x) = \left\{
\begin{array}{cc}
p(G,x)\prod_{s=1}^{(n-1)/2} \left(p^2(G,x) - 4x^2p(G_v,x)p(G_w,x) \cos^2\left(\frac{s\pi}{n+1}\right) \right) & \mbox{if $n$ is odd},\\
\prod_{s=1}^{n/2} \left(p^2(G,x) - 4x^2p(G_v,x)p(G_w,x) \cos^2\left(\frac{s\pi}{n+1}\right) \right) & \mbox{if $n$ is even},
\end{array}
\right.
\end{equation}
as long as $p^2(G,x) -4x^2p(G_v,x)p(G_w,x) > 0$.

It is well-known (see e.g. \cite{StanleyLC}) that if $f(x)$ and $g(x)$ are LC$^+$ then so is $f(x)g(x)$. This together with the observation that $\cos^2(s\pi/(n+1)) \in [0,1]$ allows us to conclude from the hypotheses of Theorem \ref{thm-tech} that $p_n(x)$ is LC$^+$ for all $n$.
 
\section{Proof of Theorem \ref{maintheorem}} \label{sec-app}

We begin with item \ref{main-thm-item-1}. Recall that $T=S_{\ell_1,\ell_2}$ is a double star consisting of adjacent vertices $v$ and $w$ with $v$ having $\ell_1$ neighbors (other than $w$), all pendant edges, and $w$ having $\ell_2$ neighbors (other than $v$), all pendant edges, $\ell_1, \ell_2 \geq 0$ arbitrary. In what follows we assume (without loss of generality) that $\ell_1 \leq \ell_2$, and we parametrize via $\ell_1=s$, $\ell_2=s+e$ with  $s, e \geq 0$ arbitrary. To avoid cluttering the notation we use $T$ for $S_{s, s+e}$.

We have $p(T,x)=(1+x)^{2s+e} +x(1+x)^{s+e} + x(1+x)^s$, $p(T_v,x)=(1+x)^{s+e}$ and $p(T_w,x)=(1+x)^s$. We begin by arguing that $p(T,x)$ is LC$^+$. We will use (repeatedly) a combination of the following ingredients:
\begin{itemize}
\item the basic fact that $(1+x)^n$ is LC$^+$ for all $n \geq 0$, 
\item the closure of the set of polynomials that are LC$^+$ under multiplication (see e.g. \cite{StanleyLC}), and 
\item the following elementary proposition, observing that verifying the LC$^+$-ness of a perturbation of an LC$^+$ polynomial requires only checking the log-concavity relations near the coefficients that have been perturbed.
\end{itemize}
\begin{prop} \label{prop-lc}
Suppose that the polynomial $f(x)=\sum_{i=0}^m b_ix^i$ is LC$^+$, that $a$ is an integer, $0 \leq a \leq m+1$, and that $A$ is a real number. If $A$ is positive then $f(x)+Ax^a$ is LC$^+$ if both $b_{a+1}^2 \geq b_{a+2}(b_a+A)$ and $b_{a-1}^2 \geq (b_a+A)b_{a-2}$ hold (here and later $b_i=0$ if $i \not \in \{0,\ldots, m\}$, so some of these relations may hold vacuously). If $A$ is negative then $f(x)+Ax^a$ is LC$^+$ if $b_a+A > 0$ and $(b_a + A)^2 \geq b_{a+1}b_{a-1}$ holds.    
\end{prop}

An alternate expression for $p(T,x)$ is
\begin{equation} \label{eq-example}
p(T,x)=(1+x)^s\left((1+x)^e\left((1+x)^s + x\right)+x\right).
\end{equation}
From the previous discussion we see that to establish that $p(T,x)$ is LC$^+$ it suffices to show that
$$
\binom{s}{2}^2 \geq \binom{s}{3}(s+1)
$$
for integer $s \geq 0$
(this allows us to deduce, via Proposition \ref{prop-lc}, that $(1+x)^s + x$ is LC$^+$, and so, via closure under products, that $(1+x)^e((1+x)^s + x)$ is LC$^+$) and that
$$
\left(\binom{s+e}{2} + e\right)^2 \geq \left(\binom{s+e}{3} + \binom{e}{2}\right)(s+e+2)
$$
for integer $s,e \geq 0$
(this allows us to deduce that $(1+x)^e((1+x)^s + x) + x = (1+x)^{s+e} + x(1+x)^e + x$ is LC$^+$, so that so also is $p(T,x)$). The first of these log-concavity relations is straightforward to verify by hand. The second requires checking that a certain two-variable polynomial of degree $4$ with $13$ monomials is non-negative at all integer points in the first quadrant; this we verify via a {\tt Mathematica} computation. Note that since we only ever perturb by linear terms, Proposition \ref{prop-lc} can be applied for all choices of $s, e \geq 0$.

We now turn to $p^2(T,x)-4qx^2p(T_v,x)p(T_w,x)$, which for convenience we denote by $f_q(x)$, and we distinguish the cases $s \geq e \geq 0$ and $e > s \geq 0$. We consider first $s \geq e \geq 0$. In this case, writing $f_q(x)$ in decreasing order of powers of $y:=(1+x)$ (to facilitate a nested presentation), we have
$$
\begin{array}{rcl}
f_q(x) & = &
y^{4s+2e} + 2xy^{3s+2e} + 2xy^{3s+e} + 
x^2y^{2s+2e} + 2x^2(1-2q)y^{2s+e} +x^2y^{2s} \\
& = &
y^{2s}(y^{e}(y^{e}(y^{s-e}(y^{e}(y^{s} +2x)+2x)+x^2 )+ 2x^2(1-2q))+x^2).
\end{array}
$$  

Using the same strategy as for $p(T,x)$, working out from the inside of this nested expression, we see that the LC$^+$-ness of $f_q(x)$ follows from the validity of each of the following relations for integers $s \geq e \geq 0$ and for $q$ as specified: first
$$
\binom{s}{2}^2 \geq \binom{s}{3}(s+2),
$$
which shows that $(1+x)^{s} +2x$ is LC$^+$; then
$$
\left(\binom{s+e}{2}+2e\right)^2 \geq \left(\binom{s+e}{3}+2\binom{e}{2}\right)(s+e+4),
$$
which shows that $(1+x)^{e}((1+x)^{s} +2x)+2x$ is LC$^+$; and then both
$$
(2s+4)^2 \geq \binom{2s}{2}+4s-2e + 1
$$
and
$$
\left(\binom{2s}{3}+2\binom{s}{2}+2\binom{s-e}{2}\right)^2 \geq \left(\binom{2s}{4}+2\binom{s}{3}+2\binom{s-e}{3}\right)\left(\binom{2s}{2}+4s-2e+1\right),
$$
which together show that $(1+x)^{s-e}((1+x)^{e}((1+x)^{s} +2x)+2x) + x^2$ is LC$^+$. 

The next step depends on the value of $q$. For $q \in [0,1/2]$ ($2x^2(1-2q) \geq 0$), we check both
$$
(2s+e +4)^2 \geq \binom{2s+e}{2} + 4s+2e + 3-4q 
$$
and
$$
\begin{array}{c}
 \left(\binom{2s+e}{3} + 2\binom{s+e}{2} + 2 \binom{s}{2} + e\right)^2  \\
 \geq  \\
 \left(\binom{2s+e}{4} + 2\binom{s+e}{3} + 2 \binom{s}{3} + \binom{e}{2}\right)\left(\binom{2s+e}{2} + 4s+2e + 3-4q\right),
\end{array}
$$
while for $q \in (1/2,1]$ ($2x^2(1-2q) < 0$) we check 
$$
\left(\binom{2s+e}{2} + 4s+2e + 3-4q\right)^2 \geq \left(\binom{2s+e}{3} + 2\binom{s+e}{2} + 2 \binom{s}{2} + e\right)(2s+e +4),
$$
and we also check that the quadratic term of $(1+x)^{2s+e}+2x(1+x)^{s+e} +2x(1+x)^{s} + x^2(1+x)^{e} + 2x^2(1-2q)$, namely $\binom{2s+e}{2} + 4s+2e + 3-4q$, is positive. This last is immediate for $s \geq 1$. If $s=0$ (and so $e=0$ since $s \geq e$) it fails for all $q \geq 3/4$. However, in the case $(s,e)=(0,0)$ we have directly that $f_q(x) = 1 + 4x + 4(1-q)x^2$ which is evidently LC$^+$ for all $q \in [0,1]$.

All this establishes that $(1+x)^{2s+e}+2x(1+x)^{s+e} +2x(1+x)^{s} + x^2(1+x)^{e} + 2x^2(1-2q)$ is LC$^+$. Finally, we check both 
$$
(2s+2e+4)^2 \geq \binom{2s+2e}{2} + 4s+6e  + 4 -4q
$$
and
$$
\begin{array}{c}
\left(\binom{2s+2e}{3}+2\binom{s+2e}{2}+2\binom{s+e}{2} + 4e(1-q)\right)^2  \\
\geq \\
\left(\binom{2s+2e}{4}+2\binom{s+2e}{3}+2\binom{s+e}{3} + \binom{2e}{2} +2\binom{e}{2}(1-2q)\right) \left(\binom{2s+2e}{2} + 4s+6e  + 4 -4q\right) 
\end{array}
$$
for $q \in [0,1]$, which establishes that $(1+x)^{2s+2e}+2x(1+x)^{s+2e} +2x(1+x)^{s+e} + x^2(1+x)^{2e} + 2x^2(1-2q)(1+x)^{e}+x^2$ and therefore that $f_q(x)$ is LC$^+$.

All nine stated relations may be verified via {\tt Mathematica} computations, completing item \ref{main-thm-item-1} of Theorem \ref{maintheorem} in the case $s \geq e \geq 0$ ($0 \leq \ell_1 \leq \ell_2 \leq 2\ell_1$). Note that since here, and in the case we are about to examine, we only ever perturb by linear and quadratic terms, and that we perturb by a linear term before perturbing by a quadratic term, all our applications of Proposition \ref{prop-lc} are valid for all $s, e \geq 0$.

We now turn to the complementary case $0 \leq s < e$ ($0 \leq \ell_1 < \ell_2/2$). Here we have
$$
\begin{array}{rcl}
f_q(x) &
= & 
y^{4s+2e} + 2xy^{3s+2e} + x^2y^{2s+2e} +
2xy^{3s+e}  + 2x^2(1-2q)y^{2s+e} +x^2y^{2s} \\
& = &
y^{2s}(y^{e}(y^{s}(y^{e-s}(y^{s}(y^{s} +2x)+x^2)+2x )+ 2x^2(1-2q))+x^2)
\end{array}
$$  
(recall $y=(1+x)$)
and we proceed as before. When we reach the polynomial
$$
(1+x)^{2s+e} + 2x(1+x)^{s+e} + x^2(1+x)^e + 2x(1+x)^s + 2x^2(1-2q)
$$
we need to verify (among other relations)
$$
\left(\binom{2s+e}{2} + 4s+2e + 3 -4q\right)^2 \geq \left(\binom{2s+e}{3} + 2\binom{s+e}{2} + e + 2\binom{s}{2}\right)(2s+e+4)
$$
for $0 \leq s < e$ and $q \in (1/2,1]$. For $s=0$, $e=1$ this relation fails for all $q > (5-\sqrt{5})/4$. It also fails for some values of $q$ close to $1$ for $s=0$, $e=2, 3$. For all other choices of $s$ and $e$ it holds for all $q \in (1/2,1]$, and indeed the entire analysis goes through exactly as in the case $s \geq e$ to show that  $f_q(x)$ is LC$^+$ for all $q \in [0,1]$ for all pairs $(s,e)$ with $s < e$ except $(0,1)$, $(0,2)$ and $(0,3)$.

For $(0,2)$ we have 
$$
f_q(x)=1+8x+(22-4q)x^2+(26-8q)x^3+(17-4q)x^4+6x^5+x^6
$$
and for $(0,3)$ we have
$$
f_q(x)=1+10x+(37-4q)x^2+(68-12q)x^3+(78-12q)x^4+(58-4q)x^5+28x^6+8x^7+x^8
$$
both of which are easily seen to by LC$^+$ for all $q\in[0,1]$.

For $(0,1)$ we have 
$$
f_q(x) = 1+6x+(11-4q)x^2+(6-4q)x^3+x^4
$$
which is LC$^+$ only for $q \leq (11-\sqrt{21})/8 \approx .802$ (in the range $q \in [0,1]$), precluding a direct application of Theorem \ref{thm-tech} (for larger $q$ we have $(6-4q)^2 < (11-4q)$).  

Recall, however, from the proof of Theorem \ref{thm-tech} (specifically from (\ref{eq-explicit_expression})) that for this choice of $s$ and $e$ we have
\begin{equation} \label{eq-pair}
p(T^n(v,w)) = \left\{
\begin{array}{cc}
(1+3x+x^2)\prod_{s=1}^{(n-1)/2} f_{q_s}(x) & \mbox{if $n$ is odd},\\
\prod_{s=1}^{n/2} f_{q_s}(x) & \mbox{if $n$ is even},
\end{array}
\right.
\end{equation}
where $q_s = \cos^2(s\pi/(n+1))$ (note $f_1(x) = 1+6x+7x^2+2x^3+x^4 > 0$).

As previously observed, the term $f_{q_s}(x)$ is LC$^+$ only for $q_s \leq (11-\sqrt{21})/8$ or equivalently 
$$
s \geq (1/\pi)(n+1)\cos^{-1} \sqrt{(11-\sqrt{21})/8} \approx .147(n+1).
$$
However, we have that
$f_{q_{s_1}}(x) f_{q_{s_2}}(x)$ equals 
$$
\begin{array}{c}
 1+12x+(58-4(q_{s_1}+q_{s_2}))x^2+(144-28(q_{s_1}+q_{s_2}))x^3+ \\(195-68(q_{s_1}+q_{s_2})+16q_{s_1}q_{s_2})x^4 + 
(144-68(q_{s_1}+q_{s_2})+32q_{s_1}q_{s_2})x^5+ \\(58-28(q_{s_1}+q_{s_2})+16q_{s_1}q_{s_2})x^6+ (12-4(q_{s_1}+q_{s_2}))x^7+x^8.
\end{array}
$$
A straightforward but tedious calculation shows that this product is LC$^+$ if $q_{s_1}, q_{s_2}$ satisfy $q_{s_1} \in [0,1]$ and $q_{s_2} \in [0, (13-\sqrt{33})/8 \approx .907]$, or equivalently
$$
s_2 \geq (1/\pi)(n+1)\cos^{-1} \sqrt{(13-\sqrt{33})/8} \approx .098(n+1).
$$
Pairing up the multiplicands in the product(s) on the right-hand side of (\ref{eq-pair}) by pairing the term corresponding to $s=1$ with that corresponding to $s=\lfloor n/2 \rfloor$, $s=2$ with $s=\lfloor n/2 \rfloor - 1$, and so on, there will always be at least one term in each pair with $s \geq (1/\pi)(n+1)\cos^{-1} \sqrt{(13-\sqrt{33})/8}$, and so the product of the terms is LC$^+$. If $\lfloor n/2 \rfloor$ is odd then the term corresponding to $s=\lceil \lfloor n/2 \rfloor/2 \rceil$ has no partner, but since $s \geq (1/\pi)(n+1)\cos^{-1} \sqrt{(11-\sqrt{21})/8}$ in this case this term is itself LC$^+$; and if $n$ is odd then the term $1+3x+x^2$ has no partner, but evidently this is LC$^+$. So $p_n(x)$ can be realized as a product of LC$^+$ polynomials (some quadratic, some quartic and some octic), and so is LC$^+$. This finishes the verification of item \ref{main-thm-item-1} of Theorem \ref{maintheorem}.

We now turn to item \ref{main-thm-item-2}. Recall that $T=P_k$ is a path on $k$ vertices, $k \geq 2$, with $w$ a leaf and $v$ its unique neighbor. There are well-known explicit expressions for $p(T,x)$, $p(T_w,x)$ and $p(T_v,x)$. Specifically
$$
p(T,x) = \sum_{j=0}^{\lfloor (k+1)/2 \rfloor} \binom{k+1-j}{j} x^j,
$$  
$$
p(T_w,x) = \sum_{j=0}^{\lfloor (k-1)/2 \rfloor} \binom{k-1-j}{j} x^j
$$
and
$$
p(T_v,x) = \sum_{j=0}^{\lfloor (k-2)/2 \rfloor} \binom{k-2-j}{j} x^j
$$
(note that $T$, $T_v$ and $T_w$ are all paths of varying lengths).

It is easy to verify that $p(T,x)$ is LC$^+$ for all $k$, so we turn attention to
\begin{equation} \label{eq-path-ver}
f_q(x) = p^2(T,x)-4qx^2p(T_v,x)p(T_w,x) := \sum_{j=0}^{2 \lfloor (k+1)/2 \rfloor} c_jx^j.
\end{equation}
To verify that $T^n(v,w)$ has log-concave independent set sequence for all $n$ it suffices (via Theorem \ref{thm-tech}) to check that for all $j=0, \ldots, 2 \lfloor (k+1)/2 \rfloor$ and $q \in [0,1]$ we have $c_j > 0$, and that for all $j=1, \ldots, 2 \lfloor (k+1)/2 \rfloor -1$ and $q \in [0,1]$ we have $c_j^2 \geq c_{j-1}c_{j+1}$.

We have not been able to verify these conditions for all $k$; however, a {\tt Mathematica} calculation shows that they hold for all $2 \leq k \leq \currentbestpath$ except $k=3, 5$. (There is also a slight anomaly at $k=2$, $q=1$, where $c_2=0$; but in this case $f_q(x)$ is easily seen to be a linear LC$^+$ polynomial.)

The case $k=3$ has already been dealt with (it is the case $s=0$, $e=1$ of item \ref{main-thm-item-1}). We deal with the case $k=5$ similarly. We have in this case that
$$
f_q(x) =   1 + 10x + (37-4q)x^2 + (62-20q)x^3 + (46-28q)x^4 + (12-8q)x^5 + x^6
$$
which is only LC$^+$ for $q \leq (41-\sqrt{113})/32 \approx .949$. However in this case it straightforward to verify that the product polynomial $f_{q_1}(x)f_{q_2}(x)$ is LC$^+$ for all $q_1, q_2 \in [0,1]$, and so we can use a simpler argument than in the case $k=3$: any partition of the multiplicands in (\ref{eq-explicit_expression}) into pairs leads to a factorization of $p(T^n(v,w))$ into LC$^+$ factors, with the only care being needed if the product in (\ref{eq-explicit_expression}) has an odd number of multiplicands, in which case one singleton block corresponding to a term with $\cos^2 s\pi/(n+1) \leq (41-\sqrt{113})/32$ or
$$
s \geq (1/\pi)(n+1)\cos^{-1} \sqrt{(41-\sqrt{113})/32} \approx .072(n+1).
$$
is required, and this is easily achieved.

\section{Proof of Proposition \ref{prop-claw}} \label{sec-rr}

For any graphs $G$ and $H$ with $H$ having distinguished vertex $v$, recall that $G_v[H]$ is obtained from $G$ by attaching a copy of $H$ at each vertex, with $v$ the point of attachment (see Proposition \ref{prop-zhu}). From \cite[Theorem 10]{Rosenfeld} we have the following identity for the independence polynomial of $G_v[H]$:
\begin{equation} \label{factor}
p(G_v[H],x) = \prod_{i=1}^{|V(G)|} (p(H-v,x)-\lambda_i x p(H-N[v],x))
\end{equation}
where the $\lambda_i$'s are the roots of the polynomial $x^np(G,1/x)$. Taking $H$ to be the star with $k$ leaves, with $v$ the center of the star, we obtain from (\ref{factor}) that
$$
p(G_k,x) = \prod_{i=1}^{|V(G)|} ((1+x)^k-\lambda_i x).
$$
If $\lambda_i$ is real then using Proposition \ref{prop-lc} and the log-concavity of $(1+x)^k$ for all $k \geq 0$ it is straightforward to verify that $(1+x)^k - \lambda_i x$ is LC$^+$ for all large enough $k$ (depending on $\lambda_i$). If $\lambda_i=a+b\sqrt{-1}$ with $a, b$ real and $b \neq 0$ then there is $j \neq i$ with $\lambda_j=a-b\sqrt{-1}$ and we examine
$$
\begin{array}{rcl}
((1+x)^k - \lambda_i x)((1+x)^k - \lambda_j x) & = & (1+x)^{2k} -2ax(1+x)^k +(a^2+b^2)x^2 \\
& = & (1+x)^k((1+x)^k -2ax) +(a^2+b^2)x^2.
\end{array}
$$
Repeated applications of Proposition \ref{prop-lc}, log-concavity of $(1+x)^k$, and the closure of LC$^+$ polynomials under multiplication show that this polynomial is LC$^+$ for all large enough $k$ (depending on $a$ and $b$); this is very similar to the verification that $p(T,x)$ in (\ref{eq-example}) is LC$^+$. That $p(G_k,x)$ is LC$^+$ now follows from one last application of closure of LC$^+$ polynomials under multiplication.  

\section{Proof of Theorem \ref{thm-pendant}} \label{sec-pendant}

Focussing on the vertex of attachment $v$ we have
$$
p(G^1_n,x) = p(G-v,x)(1+x)^n + xp(G_v,x)
$$
and focussing on $w$ we have 
$$
p(G^2_n,x) = p(G,x)(1+x)^{n-1} + xp(G-v,x).
$$
Theorem \ref{thm-pendant} thus follows from the following result, whose proof will occupy the rest of the section.
\begin{prop} \label{prop-bin}
Let $g(x)=\sum_{i=0}^\ell g_i x^i$ and $h(x)=\sum_{j=1}^{\ell'} h_j x^j$ have positive coefficients. For all sufficiently large $n=n(g,h)$ the polynomial $g(x)(1+x)^n + h(x)$ has unimodal coefficient sequence.
\end{prop} 

We begin by considering $g(x)(1+x)^n = \sum_{k=0}^{n+\ell} a_k x^k$, where
$$
a_k = g_0\binom{n}{k} + g_1\binom{n}{k-1} + \ldots + g_\ell\binom{n}{k-\ell}.
$$ 
For $k \leq \lceil n/2 \rceil$ we have $\binom{n}{k-j} \geq \binom{n}{k-j-1}$ (using the standard properties of the binomial) and so, using the positivity of the $g_i$'s, we have $a_{k-1} \leq a_k$ for all $k \leq \lceil n/2 \rceil$. Similarly we have $a_k \geq a_{k+1}$ for all $k \geq \lfloor n/2 \rfloor + \ell$.  

To establish unimodality of the coefficient sequence of $g(x)(1+x)^n$, it remains to establish the unimodality of $(a_{\lceil n/2 \rceil}, \ldots a_{\lfloor n/2 \rfloor + \ell})$. We proceed with the analysis in the case when $n$ is even; the case $n$ odd is virtually identical. We use the estimate
$$
\frac{\binom{n}{n/2-m}}{\binom{n}{n/2}} = 1-\frac{2m^2}{n} + \Theta(n^{-2})
$$
as $n \rightarrow \infty$ (with $m$ bounded), which is straightforward to verify using elementary estimates (but see also, e.g., \cite[Section 4 \& (41)]{BuricElezovic}). Applying this with $p \in [0,\ell]$ we get
\begin{equation} \label{eq-est}
\frac{n}{2}\left(1 - \frac{a_{n/2+p}}{\binom{n}{n/2}g(1)}\right) = \frac{g_0}{g(1)} p^2 + \frac{g_1}{g(1)} (p-1)^2 + \ldots + \frac{g_p}{g(1)}0^2 + \frac{g_{p+1}}{g(1)}1^2 + \ldots + \frac{g_\ell}{g(1)}(\ell-p)^2 +\Theta(n^{-1}).
\end{equation}
Notice that, without the $\Theta(n^{-1})$ error term, the right-hand side of (\ref{eq-est}) is exactly $E((X-p)^2)$ where $X$ is the random variable that takes value $i$ with probability $g_i/g(1)$. Rewriting as $p^2 -2pE(X) +E(X^2)$ it is evident that for sufficiently large $n$, as $p$ varies from $0$ to $\ell$ the right-hand side of (\ref{eq-est}) decreases to a minimum and then increases. From this it follows that $(a_{n/2}, \ldots, a_{n/2+\ell})$ is unimodal, completing the verification that $g(x)(1+x)^n$ has unimodal coefficient sequence for large $n$. 

To deal with the addition of $h(x)$ we need to show that $a_k + h_k \geq a_{k-1} + h_{k-1}$ for all $k \leq \ell'$, for which it suffices to show
$$
g_t\binom{n}{k-t} + h_k \geq g_t\binom{n}{k-1-t} + h_{k-1} 
$$
for any $t$ with $g_t \neq 0$; this is evident for $n$ sufficiently large.

\section{Concluding comments}

Given the central role that trees have played in graph theory, it is somewhat surprising that Question \ref{quest-AEMS} remains open. It is somewhat more surprising there are some simple families of trees for which we cannot answer the question. These include
\begin{itemize}
\item the family of binary (rooted) trees, and
\item the family of caterpillars (discussed after the statement of Theorem \ref{thm-pendant}). 
\end{itemize}

Though somewhat weaker than Question \ref{quest-AEMS}, verifying the truth of the following probabilistic statement would be helpful progress: if $a=(a_0,a_1,\ldots,a_n)$ is the independent set sequence of the random uniform tree on $n$ vertices (labelled or unlabelled) then the probability that $a$ is unimodal tends to $1$ as $n$ tends to infinity.

It would be of interest to find other families of trees, besides the family of paths, whose independent set sequence has the real-roots property, that could act as ``seeds'' for Propositions \ref{prop-zhu} and \ref{prop-claw}. The family of {\em Fibonacci trees}, a collection of recursively defined rooted trees \cite{Wagner}, is a candidate family. The Fibonacci tree $F_0$ consists of a single vertex, the root. The Fibonacci tree $F_1$ consists of a single edge, with one of the leaves designated the root. For $n \geq 2$  the Fibonacci tree $F_n$ is obtained from $F_{n-1}$ and $F_{n-2}$ (on disjoint vertex sets) by adding one new vertex, designated the root, joined to the roots of $F_{n-1}$ and $F_{n-2}$. In an earlier draft of this note we conjectured that $p(F_n,x)$ has the real roots property for all $n$; Bencs \cite{Bencs} has obtained a proof of this fact. 
\begin{prob}
Characterize those trees on $n$ vertices whose independent set sequence has the real-roots property.    
\end{prob}

\subsection*{Acknowledgement}

We thank Joshua Cooper for a helpful discussion.

\end{document}